\documentclass[a4paper]{amsart}

\usepackage{hyperref}

\usepackage{amsmath}
\usepackage{amsaddr}
\usepackage{graphicx}
\usepackage{caption}
\usepackage{subcaption}
\usepackage{cite}

\usepackage{textcomp}

\usepackage{longtable}
\usepackage{caption}
\usepackage{amssymb}
\usepackage{amsmath}
\usepackage{url}
\usepackage{psfrag}
\usepackage{color}
\usepackage{etoolbox}
\usepackage[detect-all]{siunitx}
\usepackage{booktabs}
\usepackage{etoolbox}
\usepackage{mathtools}

\usepackage{layouts}
\usepackage{calc}

\usepackage[normalem]{ulem}

\robustify\itshape 
\robustify\upshape 
\robustify\mdseries 
\robustify\slshape 

\DeclarePairedDelimiter{\norm}{\lVert}{\rVert}

\newcommand{\HN}{N}
\newcommand{\ark}{a}

\newcommand{\sz}[1]{{\scriptsize#1}}
\newcommand{\fsz}[1]{{\textbf{\scriptsize#1}}}
\newcommand{\rmsz}[1]{\mbox{{\sz{\rm{#1}}}}}

\newcommand{\bb}{\overline{a}}

\newcommand{\mymat}[1]{{#1}}
\newcommand{\beq}{\begin{equation}}
\newcommand{\eeq}{\end{equation}}
\newcommand{\bec}{\begin{center}}
\newcommand{\eec}{\end{center}}

\newcommand{\plone}[2]{{{\partial} #1}/{\partial {#2}}}
\newcommand{\pltwo}[2]{{{\partial}^2 #1}/{\partial {#2}^2}}

\newcommand{\n}{\noindent}

\newcommand{\TT}{T}

\newcommand{\cmm}{\alpha}

\newcommand{\MP}{M}

\newcommand{\lesim}{\,\raisebox{-0.4ex}{$\stackrel{<}{\scriptstyle\sim}$}\,}

\newcommand{\FRKC}{\fsz{FRKC} }
\newcommand{\FRKCtwo}{\fsz{FRKC2} }

\newcommand{\RKC}{\fsz{RKC} }

\newcommand{\CVODE}{\fsz{CVODE} }

\newcommand{\ROCKtwo}{\fsz{ROCK2} }
\newcommand{\RKLtwo}{\fsz{RKL2} }
\newcommand{\ROCKf}{\fsz{ROCK4} }
\newcommand{\DUMKA}{\fsz{DUMKA} }

\newcommand{\pFRKC}{\fsz{FRKC}}
\newcommand{\pFRKCtwo}{\fsz{FRKC2}}

\newcommand{\pRKC}{\fsz{RKC}}

\newcommand{\pROCKtwo}{\fsz{ROCK2}}
\newcommand{\pRKLtwo}{\fsz{RKL2}}
\newcommand{\pROCKf}{\fsz{ROCK4}}

\newcommand{\FRKCfs}{\fsz{FRKC4s} }
\newcommand{\pFRKCfs}{\fsz{FRKC4s}}
\newcommand{\FRKCf}{\fsz{FRKC4} }

\newcommand{\pFRKCxs}{\fsz{FRKC6s}}

\renewcommand*{\Re}{\operatorname{Re}}

\newcommand{\arka}{\ark}

\newcommand{\br}{\hline}
\newcommand{\mr}{\hline}

\usepackage{listings}
\usepackage{tikz}
\usepackage{longtable}

\usepackage{pgfplots}    
\usepgfplotslibrary{external} 
\pgfplotsset{compat=1.14}

\usepackage{calc}
\newlength{\theight}
\newlength{\toffset}
\newlength{\nodegap}

\begin{document}
\title{Factorized Runge--Kutta--Chebyshev Methods}

\author{Stephen O'Sullivan}

\address{School of Mathematical Sciences, Dublin Institute of Technology, Kevin Street, Dublin 8, Ireland}

\email{stephen.osullivan@dit.ie}

\begin{abstract}
  The second-order extended stability Factorized Runge--Kutta--Chebyshev (\pFRKCtwo) class of explicit schemes for the integration of large systems of PDEs with diffusive terms is presented.  \FRKCtwo schemes are straightforward to implement through ordered sequences of forward Euler steps with complex stepsizes, and easily parallelised for large scale problems on distributed architectures.
  Preserving 7 digits for accuracy at 16 digit precision, the schemes are theoretically capable of maintaining internal stability at acceleration factors in excess of 6000 with respect to standard explicit Runge-Kutta methods. The stability domains have approximately the same extents as those of \RKC schemes, and are a third longer than those of \RKLtwo schemes. Extension of \FRKC methods to fourth-order, by both complex splitting and Butcher composition techniques, is discussed.

\n A publicly available implementation of the \FRKCtwo class of schemes may be obtained from \mbox{\url{maths.dit.ie/frkc}}

\end{abstract}

\maketitle

\section{Introduction}

Factorized Runge--Kutta--Chebyshev (\pFRKC) methods are well suited to the numerical integration of problems where diffusion limits the efficiency of standard explicit techniques. In general, such systems of PDEs may be presented as semi-discrete ordinary differential equations of the form
\beq
w'=f(t,\,w) .
\eeq
Here, the associated Jacobian is assumed to have negative eigenvalues lying close to the real axis, a good approximation for many systems of interest in astrophysical contexts.

The main use of extended stability Runge--Kutta (ESRK) methods is to fill the gap between unconditionally stable but operationally complex implicit methods, and simply implemented explicit schemes which suffer from stability constraints for stiff problems. ESRK methods are particularly useful for problems involving diffusion, where the work required by standard explicit techniques goes as the inverse square of the mesh spacing, while for extended stability methods it goes as the inverse mesh spacing. ESRK explicit schemes may be broadly divided into factorized and recursive types.

Factorized ESRK methods are particularly straightforward to implement at second-order, consisting solely of forward Euler steps. At orders above two, splitting methods or, alternatively, additional finishing stages are required for nonlinear problems. First considered by \cite{saul1960integration,guillou1960domaine,gentzsch1978one}, factorized ESRK methods fell out of common usage for some time, until revived in 1996 as SuperTimeStepping~\cite{alexiades1996super} at first-order, and later extended to second-order applications in astrophysical simulations by means of Richardson extrapolation~\cite{o2006explicit,o2007three}. \DUMKA methods exist at third- and fourth-order~\cite{medovikov1998high}.

The perennial problem with the factorized formalism has been that, when a very large number of stages is used, the internal amplification factors can easily drown numerical precision. Factorized methods have, as a result, largely taken a back seat to recursive methods which manage internal stability by mapping the three-term recurrence relations for orthogonal polynomials to the stability polynomials~\cite{vanderhouwen1980}. Recursive methods have been developed up to fourth-order~\cite{verwer1996explicit,sommeijer1998rkc,abdulle2001second,abdulle2002fourth,martin2009second,meyer2012second,meyer2014stabilized}. 

In the following, a formulation of factorized methods is presented which has high internal stability, and is more straightforward to implement than recursive methods, and demonstrates comparable efficiency.

\section{\FRKC stability polynomials}

Stability polynomials form the backbone of ESRK numerical schemes and encapsulate the linear stability properties. Linearising a system of semi-discrete ODEs
\beq
\mymat{w}'=\mymat{A}w ,
\eeq
the FRKC stability polynomial~\cite{frkcjcp} is obtained by seeking a polynomial of degree $L$ which yields a forward Euler scheme of order $N$ for linear problems through its roots via
\beq
\mymat{W}^{L}=\mymat{W}^0+\TT\sum^{L}_{l=1}\ark_l f(W^{l-1})\quad \mymat{W}^0=w^n,\, \mymat{w}^{n+1}=\mymat{W}^{L} .
\eeq

The $M$-th stability polynomial of order $N$, with $L=MN$, determined via $\mymat{w}^{n+1}=R^N_M(\TT A) \mymat{w}^{n}$, must match the first $N+1$ terms in the Taylor expansion of the evolution operator
\beq
R^N_M(\TT A)={\rm e}^{\TT A}+\mathcal{O}(\TT A)^{N+1} .
\eeq

\n Equivalently, the linear order conditions may be expressed as constraints on the derivatives of the stability polynomial evaluated at zero:
\beq
R^{\HN}_M{}^{(n)}(0)=1 , \quad n=1,\,\cdots,\,N .
\label{eqn:order}
\eeq

\n In addition, stability requires that the polynomial is bounded according to
\beq
|R^N_M(z)|\le 1 \quad \forall \quad z=\TT \lambda .
\eeq
\n The objective is to determine a closed form for the  polynomial such that the extent of the stability domain along the negative real axis $\beta$ is as great as possible.

It is shown in~\cite{frkcjcp} that the \FRKC stability polynomial of rank $N$, and degree $L$, is a sum of Chebyshev polynomials of the first kind given by
\beq
B^{\HN}_M(z) = d^{\HN}_0 + 2\sum_{k=1}^{\HN} d^{\HN}_k C_{kM}(z),
\label{eqn:stab}
\eeq
where $C_{kM}$ is the Chebyshev polynomial of the first kind of degree $kM$, and the coefficients $d_k$ are determined through the linear order conditions given by Equation~\ref{eqn:order}. The resultant scheme follows immediately from the roots of the polynomial, $\zeta_l$, with coefficients given by
\beq
\ark_l=\frac{1}{\MP^2\cmm_M}~\frac{1}{1-\zeta_l} .
\eeq

\begin{figure}[!htbp]
\centering
\resizebox{0.8\columnwidth}{!}{\input{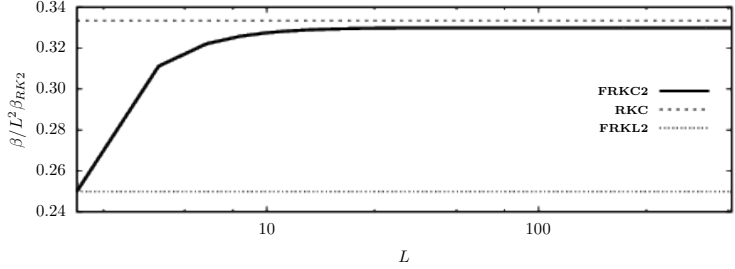}}
\caption{Extent of real stability domain at second-order for \FRKCtwo (solid line). For large values of stage-number $L$, the scheme has a stability domain which is $0.33L^2$ times the extent of the reference second-order Runge--Kutta scheme. Also shown are the corresponding values for \RKC (dashed line) and \RKLtwo (dotted line).}
\label{fig:Dvalsfig}
\end{figure}

The dependency of $\beta$ on $L$ is presented in Figure~\ref{fig:Dvalsfig} at second-order ($N=2$). While the optimal value for $\beta/L^2\beta_{\scriptsize{RK2}}$ is $0.41$~\cite{van1977}, where $\beta_{\scriptsize{RK2}}$ is the conventional second order Runge--Kutta stability limit, values of 0.330, 0.333, and 0.25 are obtained for \pFRKCtwo, \pRKC~\cite{vanderhouwen1980}, and \pRKLtwo~\cite{meyer2014stabilized} respectively.

\subsection{Damping}

\begin{figure}[!b]
  \centering
  \resizebox{0.8\columnwidth}{!}{\includegraphics[width=16cm,height=8cm]{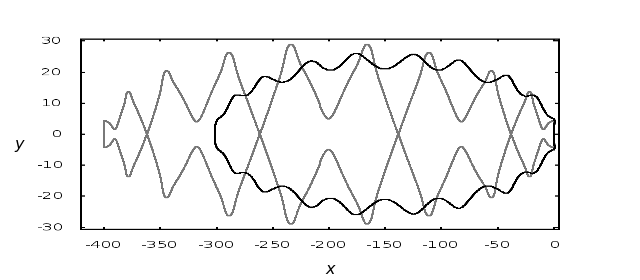}}
  \caption{Stability domain $|R(x+Iy)|=1$ for \FRKC polynomial without damping (grey line) and with sufficient damping applied for mixed hyperbolic-parabolic problems with $Pe\approx 10$ (black line). The extent of the stability domain along the real axis is contracted by approximately 25\% in the latter case.}
\label{fig:peclet}
\end{figure}

Along the real axis on the interior of the stability domain of the stability polynomial, there are points which are marginally stable, as shown in Figure~\ref{fig:peclet}. This is remedied by introducing the damping parameter $\nu\equiv\nu_0/N$ via

\beq
\ark_l=\frac{1}{(1-\nu)\MP^2\cmm_M}~\frac{1-\mu_l}{1-(1-2\mu_l)\zeta_l} ,
\eeq

\n and again enforcing the order conditions given by Equation~\ref{eqn:order} via Newton-Raphson iteration over the parameters $\mu_l$, which consist of $N$ distinct values, each repeated $M$ times. As a result, the marginally stable maxima in $|R|$ along the real axis are scaled by $\sim(1-\nu_0)$ at the expense of reducing the extent of the stability domain $\beta$ along the real axis by approximately $(1-\nu)$.  

The damping process may also be used to make the scheme applicable to problems with small hyperbolic components, with P\'{e}clet numbers $Pe\lesim 10$. For the case $Pe\approx 10$, illustrated in Figure~\ref{fig:peclet}, there is a 25\% loss in the extent of the stability domain observed.

 \subsection{Internal stability}
 Internal instability may be caused when the product of any of the possible sub-sequences of steps act to generate large values which drown out numerical precision. Following the idea of Lebedev~\cite{lebedev2000explicit,lebedev1994solve}, but with a more effective approach, the timesteps are ordered to \emph{approximately} minimise $\mathcal{Q}$, where
 \beq
 \mathcal{Q}\equiv\max(\mathcal{Q}_{j,\,k}(x)), \quad 1\le j\le  k\le L,\,\quad x\in[-\beta,\,0] , 
 \eeq
 is the maximum over the internal amplification factors defined by
 \beq
 \mathcal{Q}_{j,\,k}(x)=\prod_{l=j}^k|1+\ark_l x| .
 \eeq

\n While the optimal  value of $\mathcal{Q}$ is approximately $L^{2}$, for the purposes of constructing a stabilization algorithm, a practical upper bound of $10 L^{2}$ is chosen.  In order to maintain this bound, the estimated maximum amplification factor $\overline{\mathcal{Q}}$ is held to a minimum while  $l$ runs from 1 to $L$, where
 \beq
 \overline{\mathcal{Q}}\equiv\norm*{\max\left(\prod_{j=1}^{l}v_{j,\,k},\,\prod_{j=l+1}^{L/2}v_{j,\,k}\right)}_1 .
 \eeq

\begin{figure}[!htbp]
\centering
\resizebox{0.8\columnwidth}{!}{\input{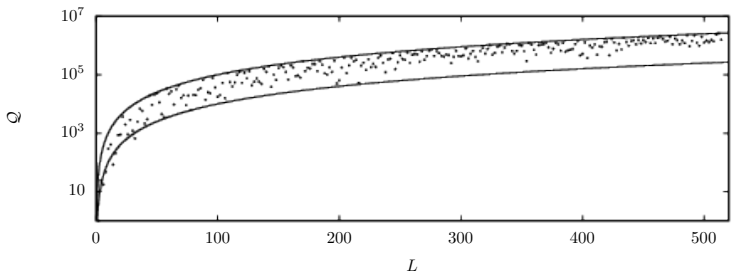}}
\caption{The maximum realised internal amplification factor $\mathcal{Q}$ as a function of $L$ for second-order schemes. Guidelines are at $L^2$ and $10L^2$.}
\label{fig:intstabfig}
\end{figure}

\n The amplification factors $v_{j,\,k}$ are defined by $v_{j,\,k}= |1+\ark_j x_k|$, where $x_k\in [\overline{\beta},\,0]$ are $L$ uniformly distributed values over the reduced range $\overline{\beta}=(1-nC)\beta$, with $C=10^{-4}$. Initially, $n=1$, however, in a limited number of cases, the process is repeated with $n$ incremented until the required bound is satisfied. In this work, the mean value of $n$ for the second-order schemes was found to be 1.5. Figure~\ref{fig:intstabfig} shows the realised values of $\mathcal{Q}$ obtained via the stabilization algorithm for second-order schemes. Preserving 7 digits for accuracy, a scheme consisting of $10^4 $ stages is therefore theoretically viable in a numerical integration carried out to 16 digit precision. 

 \subsection{Convex Monotone Property}
 
\begin{figure}[!htbp]
   \centering
   \begin{subfigure}[t]{0.5\textwidth}
     \centering
     \includegraphics[width=\linewidth]{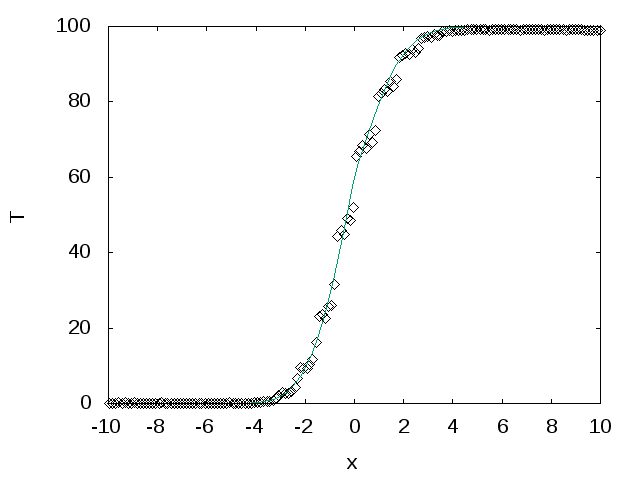}
       \caption{\RKC with 7 stages per step.}
       \label{fig:cmpa}
   \end{subfigure}
   \begin{subfigure}[t]{0.5\textwidth}
     \centering
     \includegraphics[width=\linewidth]{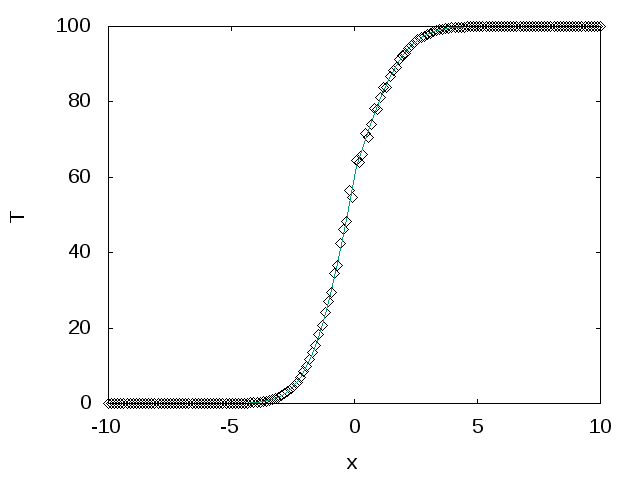}
     \caption{\FRKCtwo with 8 stages per step.}
     \label{fig:cmpc}
   \end{subfigure}
   \begin{subfigure}[t]{0.5\textwidth}
     \centering
     \includegraphics[width=\linewidth]{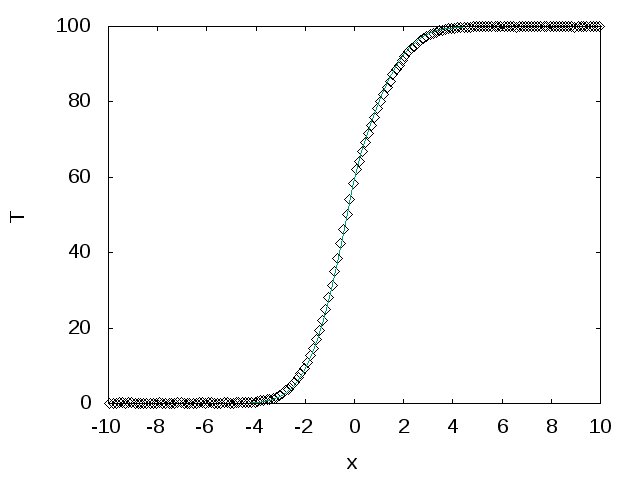}
     \caption{\RKLtwo with 7 stages per step.}
     \label{fig:cmpb}
   \end{subfigure}%
   \begin{subfigure}[t]{0.5\textwidth}
     \centering
     \includegraphics[width=\linewidth]{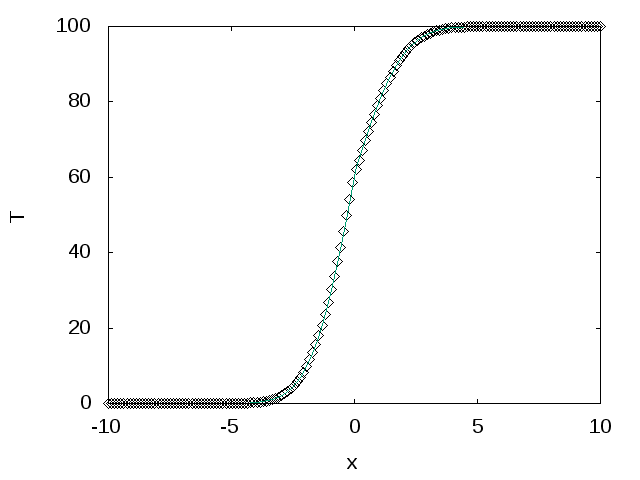}
     \caption{\FRKCtwo with error control, and same number of function evaluations as for Figure~\ref{fig:cmpc}.}
     \label{fig:cmpd}
   \end{subfigure}
   \caption{Temperature across materials of different thermal conductivities illustrating the influence of the convex monotone property (CMP)~\cite{meyer2014stabilized}.}
   \label{fig:cmp}
 \end{figure}

 The convex monotone property (CMP) is of relevance to problems with spatially varying diffusion coefficients~\cite{meyer2014stabilized}. Figure~\ref{fig:cmp} shows solutions obtained for a problem considered in~\cite{meyer2014stabilized} describing two materials placed into contact with differing temperatures and diffusion coefficients. For Chebyshev polynomial-based schemes such as \RKC and \pFRKCtwo, noise in the solutions associated with failure to meet the CMP is evident if the schemes are forced to take steps with a uniform number of stages. Figures~\ref{fig:cmpa}~and~\ref{fig:cmpc} illustrate this for the cases $L=7$ and $L=8$ respectively. The \RKLtwo scheme maintains the CMP naturally with the result that the solution is smooth (Figure~\ref{fig:cmpb}).  As shown in Figure~\ref{fig:cmpd}, for \FRKCtwo with an identical number of function evaluations as used for Figure~\ref{fig:cmpc}, but with an error-control procedure (discussed further in Section~\ref{sec:scheme}), the noise is also absent.

\section{\FRKCtwo schemes}
\label{sec:scheme}

A system $w'=\mymat{f}(w)$ is assumed such that the associated Jacobian has an eigenvalue of maximum magnitude $|\lambda|_{\rmsz{max}}$. Then, given a numerical solution $\mymat{w}^n$ at some time index $n$, $L= 2M$ stages $\mymat{W}^l$ (for $l=0,\,\cdots,\,L$) are evaluated such that $\mymat{W}^0= \mymat{w}^n$ and a second-order accurate solution $\mymat{w}^{n+1}=\mymat{W}^L$ is obtained a time $\TT$ later.  The intermediate stages of the scheme are determined via the Euler steps

\beq
\label{eqn:scheme}
\mymat{W}^{l+1} = \mymat{W}^{l}+\TT\arka_{l} \mymat{f}(\mymat{W}^{l}) . 
\eeq

Error control is straightforward since a first-order solution is available at no additional cost in function evaluations. This first-order solution $\mymat{\widehat{W}}^L$ is obtained by considering only the real parts of $\arka_l$ and $\mymat{f}(\mymat{W}^{l})$. Setting $\mymat{\widehat{W}}^{0} = \mymat{W}^{0}$, 

\beq
\label{eqn:errscheme}
\mymat{\widehat{W}}^{l+1} = \mymat{\widehat{W}}^{l}+\TT\Re({\arka_{l}}) \Re({\mymat{f}(\mymat{W}^{l})}) . 
\eeq

\n The error, scaled to a specified tolerance $TOL$,  is estimated using
\beq
\|err\|=\left\| \frac{|\mymat{W}^{l+1}-\widehat{W}^{l+1}|}{TOL(1+\max(|\mymat{W}^{l+1}|,\,|\widehat{W}^{l+1}|)}\right\| .
\eeq
\n
If $\|err\|>1$, the step is rejected and retried with $\TT$ scaled by $SAFE/\sqrt{\|err\|}$. Otherwise, a predictive controller is used to determine the subsequent timestep calibrated to the required tolerance via
\beq
\TT^{n+1}=\left(\frac{SAFE}{\sqrt{\|err^n\|}}\right)\left(\frac{\TT^n}{\TT^{n-1}}\right)\sqrt{\frac{\|err^{n-1}\|}{\|err^n\|}} ,
\eeq
using values of $\TT$ and $\|err\|$ from previous timesteps. $SAFE$ is a safety factor chosen with a value 0.8 here. The reader is referred to~\cite{hairerstiff} for further details of error control procedures.

\subsection{\FRKCtwo public code}
\label{demo}

A freely available C implementation of the second-order \FRKCtwo schemes may be accessed at \url{maths.dit.ie/frkc}. The files \href{http://www.maths.dit.ie/frkc/downloadfrkc2core.php}{frkc2core.c} and \href{http://www.maths.dit.ie/frkc/downloadfrkc2user.php}{frkc2user.c} provide the code for internal calculations required by the \FRKCtwo scheme and the code specific to the particular problem respectively. For a given value of $M$, up to 257, the extent of the stability domain along the real axis, $\overline{\beta}$, and the maximum realised internal amplification factor $\mathcal{Q}$ (see also Figure~\ref{fig:intstabfig}) are given on line  $3M-2$ of \href{http://www.maths.dit.ie/frkc/downloadfrkc2arks.php}{frkc2arks.dat}. The real and imaginary parts of $\arka_{l}$ are recorded on lines $3M-2$ and $3M$ respectively. 

\begin{figure}[!htbp]
\centering
\resizebox{0.8\columnwidth}{!}{\input{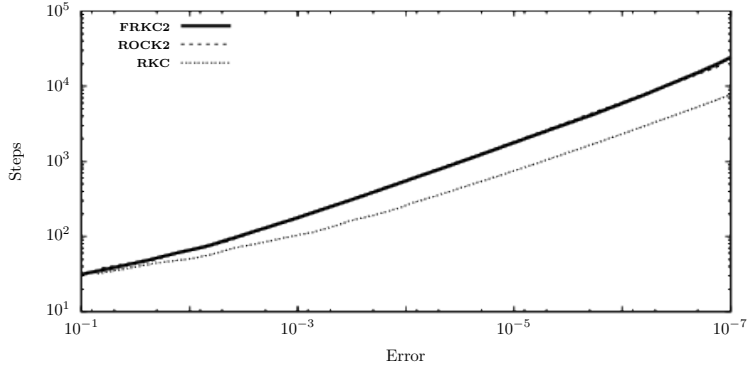}}
\caption{Efficiency comparison for the two-dimensional Brusselator problem for \pFRKCtwo, \pROCKtwo, and \pRKC. (The lines for \FRKCtwo and \ROCKtwo are almost coincident.)}
\label{fig:bruss2dcomparisonfig}
\end{figure}

The default problem provided in \href{http://www.maths.dit.ie/frkc/downloadfrkc2user.php}{frkc2user.c} is a two-species reaction diffusion Brusselator

\beq
\begin{aligned}
  \plone{v}{t}&=0.02\left(\pltwo{v}{x_1} +\pltwo{v}{x_2}\right) +1 -4v +w v^2 ,\\ 
  \plone{w}{t}&=0.02\left(\pltwo{w}{x_1} +\pltwo{w}{x_2}\right)  +3 v  -v^2w  ,\\
  v(0,\,x)&=1+\sin(2\pi x)\quad w(0,\,x)=3+\cos(2\pi y) ,
\end{aligned}
\eeq

\n  which possesses a spectral radius of approximately $6400 $ for a $200\times 200$ mesh. The initial state is a perturbation of the equilibrium solution~\cite{nicolis1977self}which is given by $v=1$, $w=3$. Figure~\ref{fig:bruss2dcomparisonfig} shows the number of steps required to attain a given error in the solution for \pFRKCtwo, \pROCKtwo~\cite{abdulle2001chebyshev}, and \pRKC. It is evident that, for a given precision, there is little difference in the number of steps required by \FRKCtwo and \pROCKtwo. At higher degrees of acceleration (fewer steps), the difference between the three schemes is negligible.

\section{\FRKCf scheme}

\subsection{Complex splitting}

Above second-order, nonlinear order conditions are present which require additional consideration. One approach, given a semi-linear parabolic (reaction-diffusion) equation of the form $w'=Aw+f_B(w)$, is to split the nonlinear part $f_B(w)$, which is typically easily integrated, from the linear diffusion terms $Aw$. For orders above two, this requires complex timesteps~\cite{Castella09,blanes2013optimized} and may be prescribed in the form

\beq
\mymat{w}^{n+1}={\rm e}^{\TT_{k_J}  \mymat{B}}{\rm e}^{\TT_{k_{J-1}}  \mymat{A}}\cdots{\rm e}^{\TT_{k_3}  \mymat{B}}{\rm e}^{\TT_{k_2}  \mymat{A}}{\rm e}^{\TT_{k_1}  \mymat{B}} \mymat{w}^n .
\eeq

\begin{figure}[!bthp]
\centering
\resizebox{0.8\columnwidth}{!}{\input{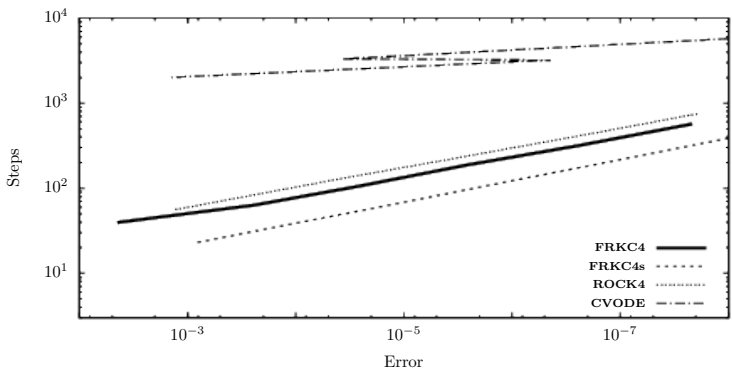}}
   \caption{Performance results derived from the estimated error for the Brusselator problem. Shown are data for the fourth-order ESRK schemes \pFRKCfs, \pFRKCxs, \pROCKf, and the fifth-order implicit \CVODE scheme.}
\label{fig:compfig}
\end{figure}

Figure~\ref{fig:compfig} shows the split \FRKCfs scheme is competitive with \pROCKf~\cite{abdulle2002fourth}. However, in support of the splitting approach, it may be noted that \ROCKf suffers significantly from internal stability issues arising from the finishing stages required for the nonlinear order conditions (discussed further in Section~\ref{sec:butcher}) which effectively limits the scheme to about 150 stages.

\subsection{Butcher composition}
\label{sec:butcher}

An alternative to complex splitting is Butcher composition. At fourth-order, as illustrated in Table~\ref{tab:butchertab}, $L-4$ forward Euler steps are adopted from the \FRKC stability polynomial with $N=4$ as the scheme $a$, and appropriate finishing stages for the scheme $\bb$ are subsequently derived. 
\begin{table}[!bp]
  \begin{minipage}{0.55\textwidth}
    \begin{center}
      \begin{tabular}{c|cccccc}
        \br
        0  &   &    &   &   &   &      \\
        $c_2$ & $a_{1}$  &    &   &   &   &      \\
        $c_3$ & $a_{1}$  & $a_{2}$   &   &   &   &      \\
        $c_4$ & $a_{1}$  & $a_{2}$   & $a_{3}$   &   &   &     \\
        $\vdots$ & $\vdots$ & $\vdots$  & $\vdots$   & $\ddots$  &   &      \\
           $c_L$ & $a_{1}$  & $a_{2}$   & $a_{3}$   & $\cdots$  & $a_{{L-5}}$  &      \\
        \mr
        & $a_1$  & $a_2$   & $a_3$   & $\cdots$  & $a_{L-5}$  & $a_{L-4}$\\
        \br
      \end{tabular}
    \end{center}
  \end{minipage}%
  \begin{minipage}{0.45\textwidth}
    \begin{center}
      \begin{tabular}{c|cccccc}
        \br
        0  &   &    &   &   \\
        $\overline{c}_2$ & $\overline{a}_{2\, 1}$  &    &   &   \\
        $\overline{c}_3$ & $\overline{a}_{3\, 1}$  & $\overline{a}_{3\, 2}$   &   &   \\
        $\overline{c}_4$ & $\overline{a}_{4\, 1}$  & $\overline{a}_{4\, 2}$   & $\overline{a}_{4\, 3}$   &   \\
        \mr
        & $\bb_1$  & $\bb_2$   & $\bb_3$   & $\bb_4$ \\
        \br
      \end{tabular}
       \end{center}
  \end{minipage}
  \caption{\label{tab:butchertab}\FRKCf order composition tableaux.}
\end{table}

According to a theorem of Hairer \& Wanner~\cite{hairernonstiff}, given the B-series $a$, $b$, the composite scheme $ab$ is determined via
\begin{equation}
  a\bb(t)=\frac{1}{\alpha(t)}\sum\left(\sum_{i=0}^q\binom{q}{i}\bb(s_i(t))\prod_{z\in d_i(t)}a(z)\right) ,
\label{eqn:comp}
\end{equation}
where rooted trees $t$ are used to represent derivatives in Taylor series. The first summation in Equation~\ref{eqn:comp} is over all $\alpha(t)$ different labelings of $t$, $s_i(t)$ is the subtree formed by the first $i$ indices, and $\{z\in d_i(t)\}$ is the difference set of subtrees formed by the remaining indices. The eight order conditions at fourth order are then given by

\beq
   \begin{aligned}
      a\bb(
  \setlength{\theight}{0ex} 
  \setlength{\toffset}{\theight-1ex} 

) .
    \end{aligned}
   \eeq
   
   \n Hence, for  given scheme $a$, imposing the required order conditions on $a\bb$ yields equations for $\bb$, which are in turn easily solved for $\bb$. The reader is referred to~\cite{hairernonstiff} for further details. 

   Figure~\ref{fig:compfig} shows a comparison of the \FRKCf scheme based on composition methods with other schemes. The reference solution is provided by a fifth-order implicit preconditioned BDF solver from the \CVODE numerical package~\cite{cohen1996cvode}. The number of steps required for a given precision is comparable for \ROCKf and \FRKCf and somewhat greater than the split \FRKCfs scheme. This difference may be ascribed to a loss of precision due to the accumulation of errors over the finishing stages~\cite{hundsdorfer2013numerical}.

   \section{Conclusions}

   \FRKC extended stability Runge--Kutta (ESRK) schemes are shown to well-suited to the integration of large-scale problems governed by systems of PDEs where the explicit timescale is constrained by diffusion. An implementation of the \FRKCtwo second-order schemes, publicly available at \url{maths.dit.ie/frkc}, is presented.
   The fourth-order \FRKCf schemes are also presented with nonlinear order conditions addressed via both splitting and composition methods. These schemes are shown to be competitive with alternative ESRK methods.

\subsection{Acknowledgments}
I am grateful to organisers of Astronum2016 for the invitation to present this work in Monterey and to an anonymous referee for helpful comments.

\section{References}
\bibliography{osullivan2016}{}
\bibliographystyle{plain}

\end{document}